\RequirePackage{ifpdf}
\ifpdf % We are running pdfTeX in pdf mode
\documentclass[pdftex]{sigma}
\else
\documentclass{sigma}
\fi

\font\frak=eufm10 scaled\magstep1
\def\goth #1{\hbox{\frak #1}}
\def\g{\goth{g}}
\def\h{\goth{h}}

\def\ad{\mathop{\mathrm{ad}}}
\def\lin{\mathop{\mathrm{lin}}}
\def\lie#1{{\mathcal{L}}_{#1}}
\def\loid{\lbrack\!\lbrack}
\def\roid{\rbrack\!\rbrack}
\def\ov{\overline}
\def\nablaC{\nabla^{\scriptscriptstyle{\mathrm{C}}}}
\def\nablaT{\nabla^{\scriptscriptstyle{\mathrm{T}}}}
\def\RT{R^{\scriptscriptstyle{\mathrm{T}}}}

\def\A{\mathcal{A}}
\def\R{\mathbb{R}}

\begin{document}

\allowdisplaybreaks

\renewcommand{\thefootnote}{$\star$}

\renewcommand{\PaperNumber}{061}

\FirstPageHeading

\ShortArticleName{Cartan Connections and Lie Algebroids}

\ArticleName{Cartan Connections and Lie Algebroids\footnote{This paper is a
contribution to the Special Issue ``\'Elie Cartan and Dif\/ferential Geometry''. The
full collection is available at
\href{http://www.emis.de/journals/SIGMA/Cartan.html}{http://www.emis.de/journals/SIGMA/Cartan.html}}}

\Author{Michael CRAMPIN}

\AuthorNameForHeading{M.~Crampin}

\Address{Department of Mathematical Physics and Astronomy, Ghent University,\\
Krijgslaan 281, B--9000 Gent, Belgium\\[1mm]
{\rm Address for correspondence:} 65 Mount Pleasant, Aspley Guise, Beds MK17~8JX, UK}
\Email{\href{mailto:Crampin@btinternet.com}{Crampin@btinternet.com}}

\ArticleDates{Received March 23, 2009, in f\/inal form June 07, 2009;  Published online June 14, 2009}

\Abstract{This paper is a study of the relationship
between two constructions associated with Cartan geometries, both of
which involve Lie algebroids:\ the Cartan algebroid, due to
[Blaom A.D., {\it Trans.\ Amer.\ Math.\ Soc.}~{\bf 358} (2006),
3651--3671], and tractor calculus [\v{C}ap~A., Gover~A.R., {\it
Trans.\ Amer.\ Math.\ Soc.}~{\bf 354} (2001), 1511--1548].}


\Classification{53B05; 53C05; 53C15}

\section{Introduction}

In 2006 Anthony D.\ Blaom published a very interesting paper \cite{Blaom}
in which he proposed a wide-ranging generalization of Cartan's idea of
a geometry as a Klein geometry deformed by curvature.  The essential
ingredient of Blaom's conception is a Lie algebroid with a linear
connection (covariant derivative) that is compatible with the
algebroid bracket in a certain sense.  Blaom calls such an algebroid a
Cartan algebroid, and~-- rather unfortunately in the circumstances~-- he calls a connection with the requisite property a Cartan
connection.  A Cartan geometry in the origi\-nal sense (perhaps one
should add `as reformulated by Sharpe \cite{Sharpe}') of a principal
bundle with a~Cartan connection form~-- what Blaom would call a
classical Cartan connection~-- gives rise to a Cartan algebroid, as one
would expect.  But there are many more examples of the latter than of
the former; for instance, a Cartan algebroid arising from a Cartan
geometry must necessarily be transitive, but Blaom gives examples of
Cartan algebroids which are not transitive.  However, in this paper I
shall concentrate on the transitive case.

In the introduction to his paper Blaom speculates, in an aside, that
there may be a relationship between his work and so-called tractor
calculi.  The concept of a tractor calculus is based on the methods
developed in the 1920s and 1930s by T.Y.~Thomas and
J.H.C.~Whitehead for the study of conformal and projective
(dif\/ferential) geometries.  These methods involve the construction in
each case of a vector bundle associated with the geometry and a
covariant derivative on its sections.  Each of these geometries is a
Cartan geometry of a special class known as parabolic.  The tractor
calculus programme is concerned with def\/ining similar structures for
all parabolic Cartan geometries.  It turns out that the relevant
vector bundles also support Lie algebroid brackets.  An encyclopaedic
account is to be found in a paper of A.~\v{C}ap and A.R.~Gover~\cite{Cap}.

One aim of this note is to spell out the exact relationship between
Blaom's theory of the Cartan algebroid of a (classical) Cartan geometry
and \v{C}ap and Gover's tractor calculus.  Of course there can be a
direct relationship, so far as Blaom's theory is concerned, only when
attention is restricted to Cartan geometries.  However, it turns out
not to be necessary to restrict one's attention further to parabolic
geometries: quite a lot of tractor calculus works for arbitrary Cartan
geometries.  Both theories involve linear connections on vector
bundles, but though one can use the same vector bundle in each case,
the two connections are dif\/ferent.  The best way of describing the
relationship between the two connections, in my view, is to regard
them as secondary constructs, which may be def\/ined in terms of a
single object canonically associated with a Cartan geometry.
This is an algebroid connection: a collection of operators with
covariant-derivative-like properties, but which involve, as it were,
dif\/ferentiation of a section of the algebroid along a section of the
algebroid. Both of the papers \cite{Blaom,Cap} discuss this algebroid
connection, but tend to regard it as dependent on the relevant linear
connection; I have reversed the dependency.

A second natural question to ask about Blaom's theory is how to
specify which Cartan algebroids derive from Cartan geometries.  Those
that do must be transitive, but there are further conditions to be
satisf\/ied, in part because algebroids derived from Cartan geometries
have more structure than just that of an algebroid with Cartan
connection (in the sense of Blaom).  This extra structure comes to the
fore in tractor calculus, so it is no real surprise that one can use
the idea of a tractor connection to prove that when suitable extra
structure is specif\/ied the Cartan algebroid does derive from a Cartan
geometry.

In Section~\ref{section2} of this note I describe a vector bundle associated with a
Cartan geometry, called the adjoint tractor bundle, on which the
events under discussion take place.  In Section~\ref{section3} I def\/ine the Lie
algebroid structure on the adjoint tractor bundle, and in Section~\ref{section4} I
specify the two linear connections, the tractor connection of \v{C}ap
and Gover and Blaom's Cartan connection.  In Section~\ref{section5} I give an
alternative derivation of the Lie algebroid structure using the Atiyah
algebroid.  In Section~\ref{section6} I prove a theorem which gives suf\/f\/icient
conditions for a Cartan algebroid to come from a Cartan geometry.  In
the f\/inal section I apply some of the ideas introduced earlier to a~discussion of Cartan space forms, which give rise to Cartan algebroids
which are symmetric in Blaom's terminology.

\looseness=-1
Most of the concepts and arguments in this note are drawn from one or
other of the two main references~\cite{Blaom,Cap}.  This seems
pretty-well inevitable given that the purpose of the exercise is to
explain how the two f\/it together. Indeed, what is new herein -- and
there are several new and I believe illuminating results to be found
below -- arises from the juxtaposition of the two apparently rather
dif\/ferent approaches of these papers, with the occasional
support of an appeal to~\cite{Sharpe}.

There are several dif\/ferent types of bracket operations to be
considered, which I want to distinguish between notationally.  I shall
reserve $[\cdot,\cdot]$ for the bracket of vector f\/ields.  I shall
denote the bracket of any of the f\/inite-dimensional Lie algebras that
occur by $\{\cdot,\cdot\}$, and use the same notation for certain
natural extensions of such brackets.  And I shall use
$\loid\cdot,\cdot\roid$ to denote a~Lie algebroid bracket.

Finally in this introduction I brief\/ly review the basic facts aboout
Cartan geometries and Lie algebroids that I need.

\subsection{Cartan geometries}\label{section1.1}

Recall that a Klein geometry is a homogeneous space $G/H$ of a Lie
group $G$, where $H$ is a closed subgroup of $G$.  I denote by $\g$
the Lie algebra of $G$ and by $\h$ the subalgebra of $\g$ which is the
Lie algebra of $H$.  A Cartan geometry on a manifold $M$,
corresponding to a Klein geometry for which $G/H$ has the same
dimension, is a principal right $H$-bundle $\pi:P \to M$ together with a~$\g$-valued 1-form $\omega$ on $P$ satisfying the following
conditions:
\begin{enumerate}\itemsep=0pt
\renewcommand{\theenumi}{$\arabic{enumi})$}
\renewcommand{\labelenumi}{\theenumi}
\item the map $\omega_p : T_pP \to \g$ is an isomorphism for
each $p \in P$;
\item $R_h^* \omega = \ad(h^{-1})\omega$ for each $h\in H$, where
$h\mapsto R_h$ is the right action; and
\item $\omega(\tilde{\xi})=\xi$ for each $\xi \in \h$, where
$\tilde{\xi}$ is the fundamental vector f\/ield on $P$ corresponding to
$\xi$.
\end{enumerate}
The form $\omega$ is the Cartan connection form.  From the second
condition it follows that for $\xi \in \h$,
$\lie{\tilde{\xi}}\omega+\{\xi,\omega\}=0$; if $H$ is connected then
the dif\/ferential version is equivalent to the original condition, and
I shall assume that $H$ is in fact connected.

The curvature $\Omega$ of a Cartan connection is the $\g$-valued
2-form on $P$ given by
\[
\Omega(X,Y)=d\omega(X,Y)+\{\omega(X),\omega(Y)\}.
\]
The value of the curvature when one of its arguments is vertical
is zero:
\[
\Omega(\tilde{\xi},Y)=
\tilde{\xi}(\omega(Y))-\omega([\tilde{\xi},Y])+\{\xi,\omega(Y)\}
=(\lie{\tilde{\xi}}\omega)(Y)+\{\xi,\omega(Y)\}=0.
\]

It is clear that in the above def\/inition $G$ has no role other than
to provide the Lie algebra $\g$. Following Sharpe \cite{Sharpe} we
may replace a Klein geometry as the model of a Cartan
geometry by the pair $(\g,\h)$ consisting of a Lie algebra $\g$ and a
subalgebra $\h$ (what Sharpe calls an inf\/initesimal Klein geometry),
together with a Lie group $H$ whose Lie algebra is $\h$, and a
representation of $H$ on $\g$ which restricts to the adjoint
representation of $H$ on $\h$. The representation of $H$ on $\g$ will
be called the adjoint representation and denoted by $\ad$. A Cartan
geometry with these data is said to be modelled on $(\g,\h)$ with
group $H$.

\subsection{Lie algebroids}\label{section1.2}

A Lie algebroid is a vector bundle $\A\to M$ equipped with a linear
bundle map $\Pi:\A\to TM$ over the identity on $M$, the anchor, and an
$\R$-bilinear skew bracket $\loid\cdot,\cdot\roid$ on $\Gamma(\A)$,
the $C^\infty(M)$-module of sections of $\A\to M$, such that
\begin{enumerate}\itemsep=0pt
\renewcommand{\theenumi}{$\arabic{enumi})$}
\renewcommand{\labelenumi}{\theenumi}
\item $\loid \sigma,f\tau\roid=f\loid\sigma,\tau\roid+\Pi(\sigma)(f)\tau$
for $\sigma,\tau\in\Gamma(\A)$ and $f\in C^\infty(M)$;
\item $\loid\cdot,\cdot\roid$ satisf\/ies the Jacobi identity;
\item as a map $\Gamma(\A)\to\Gamma(TM)$, $\Pi$ is a homomorpism:
$\Pi(\loid\sigma,\tau\roid)=[\Pi(\sigma),\Pi(\tau)]$.
\end{enumerate}
(The third of these is in fact a consequence of the other two, but is
suf\/f\/iciently important to be worth mentioning explicitly.)

A Lie algebroid is transitive if its anchor is surjective.  The subset
of $\A$ consisting of elements mapped to zero by $\Pi$ is then a
vector subbundle called the kernel of the algebroid, $\ker(\A)$.  The
restriction of the algebroid bracket to $\Gamma(\ker(\A))$ is a
Lie-algebraic bracket, that is, it is $C^\infty(M)$-bilinear (as well
as being skew and satisfying the Jacobi identity).  Indeed, a vector
bundle with a Lie-algebraic bracket on its sections is a Lie algebroid
whose anchor is the zero map.  The restriction of the Lie-algebraic
bracket to any f\/ibre of $\ker(\A)\to M$ def\/ines a Lie algebra
structure on it.  It can be shown that for a transitive algebroid $\A$
over a connected base~$M$ the Lie algebras on f\/ibres of $\ker(\A)$
over dif\/ferent points of $M$ are isomorphic, so that $\ker(\A)$ is a
Lie algebra bundle.

For a Lie algebroid $\A$, an $\A$-connection $D$ on a vector bundle
$E\to M$ is an assignment to each $\sigma\in\Gamma(\A)$ of an
operator $D_\sigma:\Gamma(E)\to\Gamma(E)$ with the following
connection-like properties:
\begin{enumerate}\itemsep=0pt
\renewcommand{\theenumi}{$\arabic{enumi})$}
\renewcommand{\labelenumi}{\theenumi}
\item[1)] $D$ is $\R$-linear in both arguments, and $C^\infty(M)$-linear
in the f\/irst;
\item[2)] $D_\sigma(f\epsilon)=fD_\sigma\epsilon+\Pi(\sigma)(f)\epsilon$
for $f\in C^\infty(M)$, $\epsilon\in\Gamma(E)$.
\end{enumerate}
Since $D_\sigma\epsilon$ is $C^\infty(M)$-linear in $\sigma$, for
$a\in\A$, $D_\sigma\epsilon(a)$ depends so far as $\sigma$ is
concerned only on~$\sigma(a)$, just as is the case for an ordinary
connection.

If $D$ is an $\A$-connection on $\A$ then $D^*$ def\/ined by
\[
D^*_\sigma\tau=D_\tau\sigma +\loid\sigma,\tau\roid
\]
is also an $\A$-connection on $\A$; it is said to be dual to $D$. The
dual of the dual is the thing one f\/irst thought of. An $\A$-connection
on $\A$ has torsion $T$ given by
\[
T(\sigma,\tau)=D_\sigma\tau-D^*_\sigma\tau
=D_\sigma\tau-D_\tau\sigma-\loid\sigma,\tau\roid.
\]
Any $\A$-connection has curvature $C$ given by
\[
C(\sigma,\tau)\epsilon=
D_\sigma D_\tau\epsilon-D_\tau D_\sigma\epsilon
-D_{\loid\sigma,\tau\roid}\epsilon.
\]
An $\A$-connection whose curvature vanishes is a homomorphism from
sections of $\A$ with the Lie algebroid bracket to operators on
$\Gamma(E)$ with the commutator bracket:
\[
D_{\loid\sigma,\tau\roid}=D_\sigma\circ D_\tau-D_\tau\circ D_\sigma.
\]
When this holds $D$ is called a representation of $\A$ on $E$.  By the
homomorphism property of~$\Pi$, for any $\sigma\in\Gamma(\A)$ the map
$s\mapsto\loid\sigma,s\roid$ maps $\Gamma(\ker(\A))$ to itself; this
evidently def\/ines a~representation of $\A$ on $\ker(\A)$, which is
called the canonical representation.  (But
$\tau\mapsto\loid\sigma,\tau\roid$ is not an $\A$-connection, so there
is no inner representation of $\A$ on itself.  Self-representations
that restrict to the canonical representation, or in other words extend
it, are therefore of some special interest.)

\section{The adjoint tractor bundle}\label{section2}

Suppose given a Cartan geometry modelled on $(\g,\h)$ with group $H$.
I denote by $\ov{\g}\to M$ the bundle associated with $P$ by the
adjoint action of $H$ on $\g$, and call it (following \cite{Cap}) the
adjoint tractor bundle of the Cartan geometry.  The adjoint tractor
bundle is to be distinguished from the adjoint bundle of the principal
bundle $P$, which is the bundle associated with $P$ by the adjoint
action of $H$ on $\h$.  The adjoint bundle, which will be denoted by
$\ov{\h}$, is a subbundle of $\ov{\g}$.

Each of these bundles is a Lie algebra bundle, with bracket derived
from the bracket of the corresponding Lie algebra $\g$ or $\h$.  The
bracket on the f\/ibres of $\ov{\g}$ will therefore be denoted by
$\{\cdot,\cdot\}$, and likewise for the f\/ibres of $\ov{\h}$ (which is
a sub-Lie-algebra bundle of $\ov{\g}$, after all).

Consider next the sections of $\ov{\g}\to M$.  Any element of
$\Gamma(\ov{\g})$ can be regarded as a $\g$-valued function on $P$
which is $H$-equivariant; that is to say, a $\g$-valued function
$\sigma$ on $P$ such that $\sigma(ph)=\ad(h^{-1})\sigma(p)$.  I shall
use the same symbol (generally $\rho$, $\sigma$ or $\tau$) for both
the section and the corresponding function; sections of $\ov{\h}$, or
sections of $\ov{\g}$ which take their values in the subbundle
$\ov{\h}$, will be denoted by $s$, $t$ and so on.  The Lie algebra
bracket of f\/ibres extends to a~Lie-algebraic bracket of sections,
$(\sigma,\tau)\mapsto \{\sigma,\tau\}$,
$\sigma,\tau\in\Gamma(\ov{\g})$.  Alternatively, thinking of sections
as $\g$-valued $H$-equivariant functions, note that the Lie algebra
bracket of two such functions is also $H$-equivariant.

It is a consequence of the equivariance condition for
$\sigma\in\Gamma(\ov{\g})$ that $\tilde{\xi}(\sigma)+\{\xi,\sigma\}=0$
for all $\xi\in\h$; when $H$ is connected this is equivalent to the
original condition.

We have at our disposal a Cartan connection form $\omega$ on $P$.  For
any vector f\/ield $X$ on $P$, $\omega(X)$ is a $\g$-valued function on
$P$; and if $X$ is $H$-invariant, then $\omega(X)$ is $H$-equivariant
and so def\/ines a section of $\ov{\g}$.  Conversely, by the isomorphism
condition, with any $\g$-valued function~$\sigma$ on $P$ we can
associate a vector f\/ield~$X$ on~$P$ such that $\omega(X)=\sigma$; and
if $\sigma$ is a section of $\ov{\g}$ then $X$ is $H$-invariant.  Thus
a Cartan connection determines a 1-1 correspondence between sections
of~$\ov{\g}$ and $H$-invariant vector f\/ields on $P$.  I shall denote
the $H$-invariant vector f\/ield corresponding to a section $\sigma$ by
$X_\sigma$:\ thus $\omega(X_\sigma)=\sigma$.

If $s\in\Gamma(\ov{\h})$, so that considered as a function
$s$ takes its values in $\h$, then $X_s$ is vertical; and
conversely.  Thus under the correspondence determined by $\omega$ the
sections of $\ov{\h}$ correspond to $H$-invariant vertical vector
f\/ields on $P$.  In fact this part of the construction is independent
of the Cartan connection: the sections of the adjoint bundle of any
principal $H$-bundle $P$ are in a~natural 1-1 correspondence with the
invariant vertical vector f\/ields on~$P$.  The correspondence may be
def\/ined as follows.  Let $V$ be a vertical $H$-invariant vector f\/ield
on $P$.  Since $V$ is vertical, for each $p\in P$ there is an element
$s(p)\in\h$ such that $V_p=\widetilde{s(p)}|_p$.  Then in
order that $V$ be $H$-invariant $s$ must be $H$-equivariant.  Of
course the Cartan connection construction reproduces this
correspondence when $s$ is a section of $\ov{\h}$ because it
satisf\/ies the condition $\omega(\tilde{\xi})=\xi$ for $\xi \in \h$.

Recall that the curvature $\Omega$ of a Cartan connection is given by
\[
\Omega(X,Y)=d\omega(X,Y)+\{\omega(X),\omega(Y)\}.
\]
If $X$ and $Y$ are $H$-invariant then the $\g$-valued function
$\Omega(X,Y)$ is easily seen to be $H$-equivariant, and so def\/ines a
section of $\ov{\g}$. We may therefore consider the curvature as a
2-form on $M$ whose value is a section of $\ov{\g}$; in this guise it
is denoted by $\kappa$. We have
\[
d\omega(X_\sigma,X_\tau)=\kappa(\pi_*X_\sigma,\pi_*X_\tau)-\{\sigma,\tau\}.
\]

\section{The Lie algebroid structure}\label{section3}

Since the bracket of $H$-invariant vector f\/ields is $H$-invariant, we
can use the vector f\/ield bracket on $P$ to def\/ine a bracket operation
$\loid\cdot,\cdot\roid$ on $\Gamma(\ov{\g})$ by
\[
[X_\sigma,X_\tau]=X_{\loid\sigma,\tau\roid},\qquad\mbox{or}\qquad
\loid\sigma,\tau\roid=\omega([X_\sigma,X_\tau]).
\]
Moreover, for any $H$-invariant vector f\/ield $X$ on $P$, $\pi_*X$ is
a well-def\/ined vector f\/ield on $M$. It is then clear that $\ov{\g}$,
equipped with the bracket $(\sigma,\tau)\mapsto\loid\sigma,\tau\roid$
and the anchor map $\Pi:\sigma\mapsto\pi_*X_\sigma$, is a Lie
algebroid.

The kernel of $\Pi$ is $\ov{\h}$.  It can be shown directly
that if $s,t\in\Gamma(\ov{\h})$ then $[X_s,X_t]=X_{-\{s,t\}}$, so that
the restriction of the Lie algebroid bracket to $\Gamma(\ov{\h})$ is
the negative of the Lie-algebraic bracket:\ $\loid s,t\roid=-\{s,t\}$.

For any $\sigma,\tau\in\Gamma(\ov{\g})$, $X_\sigma(\tau)$ is a
$\g$-valued function which is easily seen to satisfy the
$H$-equivariance condition.  The map $D$ from $\Gamma(\ov{\g})$ to
operators on $\Gamma(\ov{\g})$ given by $D_\sigma\tau=X_\sigma(\tau)$
is a~$\ov{\g}$-connection on $\ov{\g}$.  It is clear that $D$
def\/ines a derivation of the Lie-algebraic bracket:
\[
D_\rho\{\sigma,\tau\}=\{D_\rho\sigma,\tau\}+\{\sigma,D_\rho\tau\}.
\]

The curvature of $D$ vanishes:
\[
D_\rho D_\sigma\tau-D_\sigma D_\rho\tau
-D_{\loid\rho,\sigma\roid}\tau=
X_\rho (X_\sigma(\tau))-X_\sigma
(X_\rho(\tau))-[X_\rho,X_\sigma](\tau)=0.
\]
That is to say, $D$ is a representation of the Lie algebroid $\ov{\g}$
on itself. I shall call it the fundamental self-representation associated
with the Cartan connection; it is called a fundamental
$D$-operator by \v{C}ap and Gover \cite{Cap}.

On the other hand, $D$ has torsion given by
\begin{gather*}
D_\rho\sigma-D_\sigma\rho-\loid\rho,\sigma\roid
 =X_\rho(\omega(X_\sigma))-X_\sigma(\omega(X_\rho))-\omega([X_\rho,X_\sigma]) =d\omega(X_\rho,X_\sigma)\\
\phantom{D_\rho\sigma-D_\sigma\rho-\loid\rho,\sigma\roid}{} =\kappa(\Pi(\rho),\Pi(\sigma))-\{\rho,\sigma\}.
\end{gather*}
Recall that $\tilde{\xi}(\sigma)+\{\xi,\sigma\}=0$.  It follows that
if $s$ takes its values in $\h$, so that $X_s$ is vertical,
$D_s\sigma=-\{s,\sigma\}$.  Then from the torsion formula, since
$\Pi(s)=0$, we have $D_\sigma s=\loid\sigma,s\roid$; in particular,
$D_\sigma(\Gamma(\ov{\h}))\subset(\Gamma(\ov{\h}))$.  Thus
the representation $D$, restricted to $\Gamma(\ov{\h})$, reduces to
the canonical representation of $\ov{\g}$ on $\ov{\h}$ coming from the
Lie algebroid structure.  One says that $D$ is a representation
extending the canonical one.

The key properties of the adjoint tractor bundle $\ov{\g}$ are these.
\begin{enumerate}\itemsep=0pt
\item It is a Lie algebra bundle, whose f\/ibres are modelled on $\g$.
\item It is at the same time a Lie algebroid.
\item As a Lie algebroid it is transitive, with kernel $\ov{\h}$ which
is a sub-Lie-algebra bundle of $\ov{\g}$ (relative to the Lie algebra
bracket of f\/ibres) modelled on the subalgebra $\h$ of $\g$.
\item The restriction of the Lie algebroid bracket to
$\Gamma(\ov{\h})$ (which is automatically a Lie-algebraic bracket) is
the negative of the restriction of the Lie-algebraic bracket of
$\ov{\g}$.
\item It is equipped with a representation $D$ of itself on itself,
each of whose operators $D_\sigma$ maps~$\Gamma(\ov{\h})$ to itself,
and $D$ extends the canonical representation of $\ov{\g}$ on
$\ov{\h}$ to $\ov{\g}$.
\item Each operator $D_\sigma$ of the representation acts as a
derivation of the Lie-algebraic bracket, and for $s\in\Gamma(\ov{\h})$,
$D_s$ coincides with the negative of inner derivation by $s$.
\end{enumerate}

\section{Covariant derivatives}\label{section4}

In the previous section it was shown that if $s\in\Gamma(\ov{\h})$
then both $D_s\sigma+\{s,\sigma\}=0$ and \mbox{$D_\sigma
s-\loid\sigma,s\roid=0$}.  To put things another way, for
$\rho\in\Gamma(\ov{\g})$, both $D_\rho\sigma+\{\rho,\sigma\}$ and
$D_\sigma\rho-\loid\sigma,\rho\roid$ depend only on $\Pi(\rho)$.  We
may therefore def\/ine, for any vector f\/ield $U$ on $M$, operators
$\nablaT_U$ and~$\nablaC_U$ on~$\Gamma(\ov{\g})$ by
\[
\left.
\begin{array}{l}
\nablaT_U\sigma=D_\rho\sigma+\{\rho,\sigma\}\\[5pt]
\nablaC_U\sigma=D_\sigma\rho-\loid\sigma,\rho\roid
\end{array}
\right\}
\quad\mbox{for any $\rho$ such that $\Pi(\rho)=U$}.
\]
It is easy to see, from the properties of $D$, that $\nablaT$ and
$\nablaC$ are covariant derivatives; that is, by this process we
def\/ine two connections ($TM$-connections) on $\ov{\g}$.  The f\/irst of
these is the tractor connection.  The second is Blaom's Cartan
connection; to be exact, in his terminology it is the Cartan
connection on the Lie algebroid $\ov{\g}$ corresponding to the
classical Cartan connection $\omega$.

The tractor connection can be specif\/ied directly in terms of $\omega$
(that is, without reference to~$D$) as follows.  For any $u\in T_xM$ and
$\hat{u}\in T_pP$, with $\pi(p)=x$, such that $\pi_*\hat{u}=u$
\[
\nablaT_u\sigma=\hat{u}(\sigma)+\{\omega(\hat{u}),\sigma\}.
\]
It is clear that for any $\sigma\in\Gamma(\ov{\g})$ the $\g$-valued
1-form $d\sigma+\{\omega,\sigma\}$ is $H$-equivariant and vanishes on
vertical vectors, so this formulation makes sense. This is
essentially how the tractor connection is def\/ined in \cite{Cap}.

Notice that
\[
\nablaT_u\{\sigma,\tau\}
=\{\nablaT_u\sigma,\tau\}+\{\sigma,\nablaT_u\tau\}:
\]
that is, $\nablaT$ is a derivation of the Lie-algebraic bracket. The
curvature $\RT$ of $\nablaT$ is given by
\[
\RT(U,V)\sigma=\{\kappa(U,V),\sigma\}.
\]
The torsion formula leads to
\[
\nablaT_{\Pi(\sigma)}\tau-\nablaT_{\Pi(\tau)}\sigma-\loid\sigma,\tau\roid
=\{\sigma,\tau\}+\kappa(\Pi(\sigma),\Pi(\tau)).
\]
This could be read as a formula for the Lie algebroid bracket in
terms of the tractor connection and its curvature; it is so
interpreted in \cite{Cap}.

Blaom's Cartan connection can also be def\/ined directly in terms of
$\omega$, as follows.  We have
\[
D_\sigma\rho-\loid\sigma,\rho\roid
=X_\sigma(\rho)-\omega([X_\sigma,X_\rho])
=(\lie{X_\sigma}\omega)(X_\rho).
\]
This suggests that we set
\[
\nablaC_u\sigma=(\lie{X_\sigma}\omega)(\hat{u})
\]
where again $u\in T_{\pi(p)}M$ and $\hat{u}\in T_pP$ with
$\pi_*\hat{u}=u$.  It is not dif\/f\/icult to see that for any
$H$-invariant vector f\/ield $X$ on $P$, $\lie{X}\omega$ is
$H$-equivariant and vanishes on vertical vectors, so again this
makes sense.

This formulation reveals the following interesting fact about
$\nablaC$.  An $H$-invariant vector f\/ield~$X$ on~$P$ such that
$\lie{X}\omega=0$ is an inf\/initesimal automorphism of the Cartan
connection form~$\omega$ (see~\cite{Cap2}).  It follows that the
inf\/initesimal automorphisms of $\omega$ correspond precisely to the
sections of $\ov{\g}\to M$ which are parallel with respect to
$\nablaC$.  (In fact $\nablaC$ makes a brief incidental appearance
in \cite{Cap2}, which is in part a study of inf\/initesimal automorphisms
of Cartan connection forms.)

Notice also that $\nablaC_U\sigma=D^*_\rho\sigma$ where $D^*$ is the
$\ov{\g}$-connection dual to $D$.  Indeed, one could almost sum up the
dif\/ference between $\nablaT$ and $\nablaC$ by saying that the f\/irst
embodies $D$ and the second $D^*$.

The fact that $D$ is a representation translates into the following
formula for $\nablaC$:
\[
\nablaC_U(\loid\sigma,\tau\roid)=
\loid\nablaC_U\sigma,\tau\roid+\loid\sigma,\nablaC_U\tau\roid
+\nablaC_{\Pi(D_{\tau}\rho)}\sigma-\nablaC_{\Pi(D_{\sigma}\rho)}\tau
\]
where $U=\Pi(\rho)$.  In \cite{Blaom} Blaom says that a connection
$\nabla$ on a Lie algebroid $\A$ is compatible with the Lie algebroid
bracket if
\[
\nabla_U(\loid\sigma,\tau\roid)=
\loid\nabla_U\sigma,\tau\roid+\loid\sigma,\nabla_U\tau\roid
+\nabla_{\nabla^*_\tau U}\sigma-\nabla_{\nabla^*_\sigma U}\tau.
\]
Here $\nabla^*$ is the operator def\/ined by
$\nabla^*_\sigma U=\Pi(\nabla_U\sigma)+[\Pi(\sigma),U]$;
it is an $\A$-connection on~$TM$ which is in a sense dual to $\nabla$.
Evidently ${\nablaC_\sigma}^*\Pi(\rho)=\Pi(D_\sigma\rho)$.  So the
fact that $D$ is a~representation implies that $\nablaC$ is compatible
with $\loid\cdot,\cdot\roid$:\ $\nablaC$ is indeed a Cartan connection
in the sense of Blaom.

It is straightforward to calculate the torsion $T^*$ and curvature
$R^*$ of the dual $D^*$ of an $\A$-connection on $\A$ in terms of $T$
and $R$, those of the original $\A$-connection $D$; the formulas are
given in \cite{Blaom}.  One f\/inds that $T^*=-T$, and in the case of
interest, where $D$ is a representation,
\[
R^*(\rho,\sigma)\tau
=D_\tau(T(\rho,\sigma))-T(D_\tau\rho,\sigma)-T(\rho,D_\tau\sigma)
=(D_\tau T)(\rho,\sigma).
\]
Given the relation between $T$ and the curvature $\kappa$, and the
fact that $D$ def\/ines a derivation of the Lie-algebraic bracket, one
can in ef\/fect replace $T$ by $\kappa$ on the right-hand side. These
results can be used to derive further properties of $\nablaC$ if
desired.

The dif\/ference tensor relating the two connections is the curvature:
\[
\nablaT_{\Pi(\rho)}\sigma-\nablaC_{\Pi(\rho)}\sigma
=\kappa(\Pi(\rho),\Pi(\sigma)).
\]

\section{The Atiyah algebroid and the Cartan connection}\label{section5}

Any principal bundle has canonically associated with it a Lie
algebroid, its Atiyah algebroid, which in the case of a principal
$H$-bundle $P\to M$ is def\/ined as follows.  The action of $H$ on $P$
lifts in a natural way to an action on the tangent bundle $TP$.  As a
manifold the Atiyah algebroid of $P$ is $TP/H$, the quotient of $TP$
under this action.  It is a vector bundle over $P/H=M$, each of whose
f\/ibres is obtained by identifying the tangent spaces $T_pP$ at points
$p$ on the same f\/ibre of $P\to M$.  The sections of $TP/H\to M$ may be
identif\/ied with the $H$-invariant vector f\/ields on~$P$; the Lie
algebroid bracket is the ordinary bracket of such vector f\/ields; the
anchor maps an $H$-invariant vector f\/ield on $P$ to its projection
onto $M$.  The sections of the kernel of the Atiyah algebroid can
therefore be identif\/ied with the $H$-invariant vector f\/ields on $P$
which are vertical (with respect to $\pi:P\to M$).  As was pointed out
before, any vertical $H$-invariant vector f\/ield~$V$ on~$P$ takes the
form $V_p=\widetilde{s(p)}|_p$ for some $H$-equivariant $\h$-valued
function~$s$ on~$P$.  So the kernel of the Atiyah algebroid can be
identif\/ied with the adjoint bundle of $P$, that is, $\ov{\h}$, the
vector bundle associated with $P$ by the adjoint action of $H$ on
$\h$.

It is clear from the discussion in Section~\ref{section3} that as a Lie algebroid
$\ov{\g}$ is isomorphic to the Atiyah algebroid, where the isomorphism
is def\/ined by the Cartan connection form~$\omega$.  In fact $\omega$
def\/ines an isomorphism, over the identity on~$P$, of~$TP$ with
$P\times\g$, which is equivariant with respect to the actions of~$H$
on the two spaces, and therefore passes to the quotients under these
actions.  As it happens, Blaom in~\cite{Blaom} uses the Atiyah
algebroid rather than the adjoint tractor bundle in his discussion of
the association of a Lie algebroid with a Cartan geometry.

It might indeed be considered preferable to describe the constructions
in Sections~\ref{section2} and~\ref{section3} in a dif\/ferent way which takes advantage of these
observations.  Suppose given the data for a~model geometry as
specif\/ied in the Section~\ref{section1.1}, namely an inf\/initesimal Klein geometry
$(\g,\h)$, a Lie group $H$ realising $\h$, and the adjoint
representation of $H$ on $\g$.  Take a principal $H$-bundle $P\to M$
where $\dim{M}=\dim{\g}-\dim{\h}$.  In the absence of a Cartan
connection form one can then construct two vector bundles over~$M$:
the adjoint tractor bundle $\ov{\g}$ (the bundle associated with~$P$
by the adjoint action of $H$ on $\g$), and the Atiyah bundle $TP/H$.
These bundles have the same f\/ibre dimension
$\dim{\g}=\dim{M}+\dim{\h}$.  The f\/irst is a Lie algebra bundle, the
second a~Lie algebroid.  It is worth emphasising that, given the model
data, each of these vector bundles has the stated additional structure
regardless of the existence of a Cartan connection form.  A Cartan
connection form, however, def\/ines a vector bundle isomorphism between
the two, which acts as the identity on $\ov{\h}$ (which can be
regarded as a subbundle of each bundle).  This isomorphism can be used
to carry the Lie algebroid structure from the Atiyah algebroid over to
the adjoint tractor bundle -- and equally the Lie algebra bundle
structure in the opposite direction.  In fact the existence of a
Cartan connection is equivalent to the existence of a vector bundle
isomorphism, over the identity, between these two bundles, which acts
as the identity on the common subbundle $\ov{\h}$.

\begin{theorem}\label{theorem1}
With model data as above, let $P\to M$ be a principal $H$-bundle where
$\dim{M}=\dim{\g}-\dim{\h}$.  Consider two vector bundles over $M$,
with the same fibre dimension: the adjoint tractor bundle $\ov{\g}$
and the Atiyah bundle $TP/H$.  The Cartan geometries on $P$ modelled
on $(\g,\h)$ with group $H$ are in $1\text{-}1$ correspondence with the vector
bundle isomorphisms $TP/H\to\ov{\g}$, over the identity, which reduce
to the identity on $\ov{\h}$.
\end{theorem}

\begin{proof}
It remains to show that an isomorphism of the bundles def\/ines a Cartan
connection on~$P$.  Let $\hat{\omega}:TP/H\to\ov{\g}$ be a vector
bundle isomorphism over the identity.  We can identify the sections of
$TP/H\to M$ with the $H$-invariant vector f\/ields on $P$, and the
sections of $\ov{\g}\to M$ with the $H$-equivariant $\g$-valued
functions on $P$.  Then $\hat{\omega}$ associates with every
$H$-invariant vector f\/ield $X$ on $P$ an $H$-equivariant $\g$-valued
function $\omega(X)$ on $P$, such that $\omega(X)$ depends
$C^\infty(M)$-linearly on its argument.  It follows that
$\omega(X)(p)$ depends only on the value of $X$ at $p$, and does so
linearly; thus for each $p$ we have a linear map $\omega_p:T_pP\to\g$,
which is an isomorphism.  In other words, $\omega$ is a $\g$-valued
1-form on $P$.  The condition that $R_h^* \omega = \ad(h^{-1})\omega$
for each $h\in H$ follows from the fact that $\omega(X)$ is
$H$-equivariant whenever $X$ is $H$-invariant.  The assumed property
of $\hat{\omega}$ that it reduces to the identity on $\ov{\h}$ is
equivalent to the condition that $\omega(\tilde{\xi})=\xi$ for
$\xi\in\h$.
\end{proof}

These observations suggest an alternative way of def\/ining a Cartan
connection form, in ef\/fect as an isomorphism of vector bundles.  To
suggest that one might def\/ine a Cartan geometry in this way may seem
to be a somewhat idiosyncratic procedure.  It is done to emphasise the
fact that the approaches of Blaom on the one hand, and \v{C}ap and
Gover on the other, concentrate on dif\/ferent aspects of the full
story.

\section{Cartan algebroids and Cartan geometries}\label{section6}

I listed earlier the properties of the adjoint tractor bundle
associated with a Cartan geometry.  I now show that for a certain
class of pairs of Lie algebras $(\g,\h)$ these properties are
def\/initive of Cartan geometries.

\begin{theorem}\label{theorem2}
Let $\g$ be a Lie algebra and $\h$ a subalgebra of $\g$ such that
\begin{enumerate}\itemsep=0pt
\renewcommand{\theenumi}{$\arabic{enumi})$}
\renewcommand{\labelenumi}{\theenumi}
\item all derivations of $\g$ are inner;
\item the centre of $\g$ is $\{0\}$;
\item the normalizer of $\h$ in $\g$ is $\h$.
\end{enumerate}
Let $H$ be the subgroup of the Lie group of inner automorphisms of
$\g$ which leaves $\h$ invariant; it is a Lie group with Lie algebra
$\h$.  Let $\ov{\g}\to M$ be a transitive Lie algebroid equipped with
a~self-representation $D$ which extends the canonical representation,
such that
\begin{enumerate}\itemsep=0pt
\renewcommand{\theenumi}{$\arabic{enumi})$}
\renewcommand{\labelenumi}{\theenumi}
\item $\ov{\g}$ is at the same time a Lie algebra bundle modelled on $\g$;
\item $\ov{\h}$, the kernel of $\ov{\g}$, is a sub-Lie-algebra bundle
of $\ov{\g}$ modelled on the subalgebra $\h$ of $\g$;
\item the restriction of the Lie algebroid bracket to $\Gamma(\ov{\h})$
is the negative of the restriction of the Lie-algebraic bracket;
\item each operator $D_\sigma$ of the representation acts as a
derivation of the Lie-algebraic bracket, and for $s\in\Gamma(\ov{\h})$,
$D_s$ coincides with the negative of inner derivation by $s$.
\end{enumerate}
Then there is a principal $H$-bundle $\pi:P\to M$ and a Cartan
connection form $\omega$ on $P$ such that~$\ov{\g}$ is the adjoint tractor
bundle of $P$, the Lie algebroid bracket of $\ov{\g}$ coincides with
that determined by $\omega$, and $D$ is the fundamental
self-representation associated with $\omega$.
\end{theorem}

The proof of this theorem is closely modelled on the proofs of
Proposition~2.3 and Theorem~2.7 of \cite{Cap}.

\begin{proof}
Let $P$ be the set of all pairs $(x,p)$ where $x\in M$ and $p$ is a
Lie algebra isomorphism $\g\to\ov{\g}_x$ such that $p(\h)=\ov{\h}_x$.
Then $P$ can be given a manifold structure such that $(x,p)\mapsto x$
is a smooth submersion $\pi$, and $\pi:P\to M$ admits local smooth
sections, so it is a f\/ibre bundle over $M$.  For each $h\in H$,
$(x,p)\mapsto (x,p\circ\ad h)$, where $\ad$ denotes the action of $H$
on $\g$, def\/ines a right action of $H$ on $P$ which makes $P$ into a
principal $H$-bundle.  The map $P\times\g\to\ov{\g}$ by
$(x,p,\xi)\mapsto p(\xi)\in\ov{\g}_x$ induces an identif\/ication of the
adjoint tractor bundle of $P$ with $\ov{\g}$.

The sections of $\ov{\g}\to M$ may be identif\/ied with the
$H$-equivariant $\g$-valued functions on $P$, as before.  It follows
from the properties of $D$ that $D_\rho\sigma+\{\rho,\sigma\}$ depends
only on $\Pi(\rho)$.  We may therefore def\/ine a covariant derivative
operator $\nablaT$ on $\Gamma(\ov{\g})$ by
\[
\nablaT_U\sigma=D_\rho\sigma+\{\rho,\sigma\}
\quad \mbox{for any $\rho$ such that $\Pi(\rho)=U$.}
\]
Then $\nablaT_U$ is a derivation of the Lie-algebraic bracket on
$\Gamma(\ov{\g})$.  Now for any $v\in T_pP$ and any $\sigma\in
\Gamma(\ov{\g})$, $\nablaT_{\pi_*v}\sigma-v(\sigma)$ is an element of
$\g$ which depends $\R$-linearly on $v$ and on $\sigma(p)$.  So we may
def\/ine an $\R$-linear map $\Phi_p:T_pP\to\lin(\g,\g)$ (the space of
$\R$-linear maps from $\g$ to $\g$) by
$\Phi_p(v)(\xi)=\nablaT_{\pi_*v}\sigma-v(\sigma)$ for any
$\sigma\in\Gamma(\ov{\g})$ such that $\sigma(p)=\xi$.  Moreover, for
each $v\in T_pP$, $\Phi_p(v)$ is a derivation of $\g$, and by
assumption~(1) about $(\g,\h)$ we may write
$\Phi_p(v)(\xi)=\{\omega_p(v),\xi\}$.  Then $\omega_p:T_pP\to\g$ is an
$\R$-linear map, as a consequence of assumption~(2) about $(\g,\h)$.

Take $v=\tilde{\xi}_p$ for any $\xi\in\h$. Then $\pi_*v=0$, and
$v(\sigma)=\tilde{\xi}(\sigma)(p)=-\{\xi,\sigma(p)\}$. Thus
$\{\omega_p(\tilde{\xi}),\sigma(p)\}=\{\xi,\sigma(p)\}$, so that
$\omega_p(\tilde{\xi})=\xi$, again by assumption~(2).

For any $H$-invariant vector f\/ield $X$ on $P$ and any
$\sigma\in\Gamma(\ov{\g})$, $X(\sigma)\in\Gamma(\ov{\g})$ also.  Thus
$\nablaT_{\pi_*X}\sigma-X(\sigma)=\{\omega(X),\sigma\}$ where
$\omega(X)\in\Gamma(\ov{\g})$.  So $\omega$ is a smooth $\g$-valued
1-form on $P$ such that $\omega(X)$ is $H$-equivariant whenever $X$ is
$H$-invariant, whence $R_h^*\omega=\ad(h^{-1})\omega$ for all $h\in H$
(again, appealing implicitly to assumption~(2)).

Suppose that $\omega_p(v)=0$ for some $v\in T_pP$.  Let
$\rho\in\Gamma(\ov{\g})$ be such that $\Pi(\rho)(x)=\pi_*v$,
$x=\pi(p)$.  Then
\[
\nablaT_{\pi_*v}\sigma-v(\sigma)=
D_\rho\sigma(p)+\{\rho(p),\sigma(p)\}-v(\sigma)=0
\]
for all $\sigma\in\Gamma(\ov{\g})$.  In particular, this holds for
$\sigma=s\in\Gamma(\ov{\h})$.  Now by assumption $D_\rho s(p)\in\h$,
and evidently $v(s)\in\h$; thus $\{\rho(p),s(p)\}\in\h$,
and so $\{\rho(p),\h\}\subset\h$.  It follows from assumption~(3)
that $\rho(p)\in\h$.  But then $\pi_*v=\Pi(\rho)(x)=0$, so $v$ is
vertical; and we know that $\omega_p$ is an isomorphism of the
vertical subspace of $T_pP$ with $\h$.  So $v=0$.  It follows that for
each $p\in P$, $\omega_p:T_pP\to\g$ is an injective map between vector
spaces of the same dimension and is therefore an isomorphism.

Thus $\omega$ is a Cartan connection form. For each
$\sigma\in\Gamma(\ov{\g})$ there is a unique $H$-invariant vector f\/ield
$X_\sigma$ on $P$ such that $\omega(X_\sigma)=\sigma$. Then
\[
\nablaT_{\pi_*X_\sigma}\tau=X_\sigma(\tau)+\{\sigma,\tau\}
=\nablaT_{\Pi(\sigma)}\tau=D_\sigma\tau+\{\sigma,\tau\},
\]
so $D_\sigma\tau=X_\sigma(\tau)$. Since $D$ is a representation we
have
\[
D_{\loid\rho,\sigma\roid}\tau=X_{\loid\rho,\sigma\roid}(\tau)
=[X_\rho,X_\sigma](\tau).
\]
Thus
\[
\loid\rho,\sigma\roid=\omega(X_{\loid\rho,\sigma\roid})
=\omega([X_\rho,X_\sigma]);
\]
but the latter expression is just the Lie algebroid bracket of $\rho$
and $\sigma$ def\/ined by the Cartan connection form $\omega$.
\end{proof}

Any Lie algebroid equipped with a Cartan connection, that is, a
connection which is compatible with the Lie algebroid bracket in the
sense given in Section~\ref{section4}, is a Cartan algebroid, in Blaom's
terminology; Cartan algebroids are to be regarded as generalizations,
at an inf\/initesimal level, of Cartan geometries.  Each Cartan geometry
gives rise to a transitive Cartan algebroid, by the construction
described in Section~\ref{section4}, but conditions for a Cartan algebroid to come
from a Cartan geometry are not given by Blaom.  Clearly, one of those
conditions must be that the algebroid is transitive.  Now a Cartan
connection $\nablaC$ on any Lie algebroid $\A$ determines an
$\A$-connection $D$: we have merely to reverse the def\/inition of
$\nablaC$ given earlier, and set
$D_\sigma\rho=\nablaC_{\Pi(\rho)}\sigma+\loid\sigma,\rho\roid$.  Then
$D$ is a self-representation such that $D_\sigma s=\loid\sigma,s\roid$
when $\Pi(s)=0$.  Thus for a transitive algebroid the existence of a
Cartan connection is equivalent to the existence of a
self-representation extending the canonical representation, as Blaom
points out.  That is to say, a transitive Cartan algebroid is simply a
transitive Lie algebroid with a self-representation extending the
canonical representation.  The theorem above therefore provides an
answer to the question of which Cartan algebroids correspond to Cartan
geometries.

\section{Symmetric Cartan algebroids}\label{section7}

There is one situation in which it is possible to derive a
Lie-algebraic structure from a Cartan algebroid:\ this occurs when the
latter is what Blaom calls (locally) symmetric.  A Cartan algebroid is
symmetric when its Cartan connection $\nablaC$ is f\/lat.  Such an
algebroid is, locally at least, isomorphic to an action algebroid,
that is, to an algebroid generated by an inf\/initesimal action of a Lie
algebra on a manifold.  The term `symmetric' comes from this latter
property, when the algebra is conceived as an algebra of inf\/initesimal
symmetries.

When the algebroid is transitive one may work with the
self-representation $D$ rather than the Cartan connection,
and in view of remarks made in Section~\ref{section4} the condition for symmetry
is that the torsion $T$ of $D$ should satisfy $DT=0$.

\begin{theorem}\label{theorem3}
Let $\A\to M$ be a transitive Cartan algebroid, with
self-representation $D$, which is locally symmetric.  Set
$\{\sigma,\tau\}=-T(\sigma,\tau)$ where $T$ is the torsion of $D$:
then $\{\cdot,\cdot\}$ is a Lie-algebraic bracket on $\Gamma(A)$ such
that
\begin{enumerate}\itemsep=0pt
\renewcommand{\theenumi}{$\arabic{enumi})$}
\renewcommand{\labelenumi}{\theenumi}
\item $\A$ is a Lie algebra bundle, whose Lie-algebraic bracket is
$\{\cdot,\cdot\}$;
\item the restriction of the Lie algebroid bracket to $\Gamma(\ker(\A))$
is the negative of the restriction of $\{\cdot,\cdot\}$;
\item each operator $D_\sigma$ of the self-representation acts as a
derivation of $\{\cdot,\cdot\}$, and for $s\in\Gamma(\ker(\A))$,
$D_s$ coincides with the negative of inner derivation by $s$.
\end{enumerate}
\end{theorem}

\begin{proof}
Evidently $-T(\sigma,\tau)$ is $C^\infty(M)$-bilinear and skew.  Any
$\A$-connection $D$ satisf\/ies a~Bianchi identity, which when its
curvature vanishes is just a dif\/ferential condition on its torsion,
namely
\[
\bigoplus_{\rho,\sigma,\tau}\left(D_\rho
T(\sigma,\tau)+T(T(\rho,\sigma),\tau)\right)=0,
\]
where $\oplus$ stands for the cyclic sum. When $DT=0$ therefore,
$-T$, thought of as def\/ining a~bracket, satisf\/ies the Jacobi
identity, and so $\{\cdot,\cdot\}$ is a Lie-algebraic bracket on
$\Gamma(A)$.

Consider (following \cite{Blaom}) those sections of $\A$ which are
parallel with respect to the Cartan connection $\nablaC$.  By the
properties of f\/lat connections the set of such sections, say $\A_0$,
forms a~f\/inite-dimensional real vector space whose dimension is the
f\/ibre dimension of $\A$.  Elements of $\A_0$ are sections $\sigma$ for
which $D_\sigma\tau=\loid\sigma,\tau\roid$ for every $\tau$ (even
those $\tau\in\Gamma(\ker(\A))$, since by assumption~$D$ extends the
canonical representation).  It follows from the fact that $D$ is a~representation that~$\A_0$ is closed under the algebroid bracket.  The
restriction to~$\A_0$ of $-\loid\cdot,\cdot\roid$ gives it the
structure of a Lie algebra.  Moreover, each f\/ibre of $\A\to M$
acquires the same Lie algebra structure: the bracket of $a,b\in\A_x$
is $-\loid\alpha,\beta\roid(x)$ where $\alpha,\beta$ are the unique
elements of~$\A_0$ such that $\alpha(x)=a$, $\beta(x)=b$.  Thus $\A$
is a Lie algebra bundle modelled on the Lie algebra~$\A_0$.  For any
$\alpha,\beta\in\A_0$,
\begin{gather*}
T(\alpha,\beta)=D_\alpha\beta-D_\beta\alpha-\loid \alpha,\beta\roid
=\loid \alpha,\beta\roid-\loid \beta,\alpha\roid-\loid \alpha,\beta\roid
=\loid \alpha,\beta\roid,
\end{gather*}
so that the restriction of the algebroid bracket to $\A_0$ coincides
with the restriction of~$T$. Take a (local) basis $\{\alpha_i\}$ of
$\A_0$ (over $\R$).  By dimension, $\A_0$ spans $\Gamma(\A)$ over
$C^\infty(M)$; that is, every $\sigma\in\Gamma(\A)$ may be uniquely
written in the form $\sigma^i\alpha_i$ (summation over $i$ intended)
with $\sigma^i\in C^\infty(M)$.  Then
\[
T(\sigma,\tau)=\sigma^i\tau^jT(\alpha_i,\alpha_j)
=\sigma^i\tau^j\loid\alpha_i,\alpha_j\roid.
\]
Thus the Lie-algebraic bracket on $\A$ arising from the fact that it
is a Lie algebra bundle is $\{\cdot,\cdot\}$.

Now take $s,t\in\Gamma(\ker(\A))$.  Then since $\Pi(s)=\Pi(t)=0$,
\[
\loid s,t\roid=s^it^j\loid\alpha_i,\alpha_j\roid=T(s,t).
\]
Thus on $\Gamma(\ker(\A))$ the algebroid bracket coincides with $T$,
and therefore with $-\{\cdot,\cdot\}$.

The condition $DT=0$ becomes
\[
D_\rho\{\sigma,\tau\}=\{D_\rho\sigma,\tau\}+\{\sigma,D_\rho\tau\}.
\]
Since $D$ extends the canonical representation of $\A$ on $\ker(\A)$,
for $s\in\Gamma(\ker(\A))$, $D_\sigma s=\loid\sigma,s\roid$.  Thus
\[
T(s,\sigma)=D_s\sigma-D_\sigma s-\loid s,\sigma\roid=D_s\sigma.
\]
That is to say, $D_s\sigma=-\{s,\sigma\}$ for all $s\in\Gamma(\ker(\A))$.
\end{proof}

In view of the correspondence mentioned earlier between sections of
$\A$ which are parallel with respect to $\nablaC$ and inf\/initesimal
automorphisms of the Cartan connection form $\omega$, we see that a~transitive Cartan algebroid which is locally symmetric in Blaom's
sense is one in which the Lie algebra of inf\/initesimal automorphisms
of $\omega$ (which is $\A_0$) has maximal dimension.  Since $\A_0$ is
closed under the Lie algebroid bracket we may write
$\loid\alpha_i,\alpha_j\roid=C_{ij}^k\alpha_k$; the $C^k_{ij}$ are the
structure constants of $\A_0$ relative to the basis $\{\alpha_i\}$.
The anchor map, which restricts to a homomorphism of $\A_0$ into
$\Gamma(TM)$, def\/ines an inf\/initesimal action of $\A_0$ on $M$.  We
have
\[
\loid\sigma,\tau\roid =
\left(\Pi(\sigma)(\tau^k)-\Pi(\tau)(\sigma^k)
+\sigma^i\tau^jC_{ij}^k\right)\alpha_k,\qquad
\sigma=\sigma^k\alpha_k,\qquad \tau=\tau^k\alpha_k.
\]
This realises $\A$ as an action algebroid corresponding to the Lie
algebra $\A_0$.

Recall that if one starts with a Cartan geometry one f\/inds that the
torsion of the fundamental self-representation is given by
\[
T(\sigma,\tau)=\kappa(\Pi(\sigma),\Pi(\tau))-\{\sigma,\tau\}.
\]
A f\/lat Cartan geometry, that is, one for which $\kappa=0$, is of
course symmetric, and the Lie-algebraic bracket constructed by the
procedure in the theorem coincides with the original one.  But $DT$
can be zero without $\kappa$ being zero, and then the new
Lie-algebraic bracket is dif\/ferent from the original one.  In such a
case we f\/ind ourselves dealing with what Sharpe in \cite{Sharpe} calls
Cartan space forms and mutation.

The curvature $\Omega$ of a Cartan geometry is constant, according to
Sharpe, if $\Omega(X,Y)$ (a $\g$-valued function on the principal
$H$-bundle $P$) is constant for every pair of vector f\/ields $X$, $Y$
on~$P$ such that $\omega(X)$ and $\omega(Y)$ are constant.  A Cartan
geometry whose curvature is constant is a~Cartan space form.

Until now I have for the most part restricted my attention to vector
f\/ields on $P$ which are $H$-invariant, and vector f\/ields $X$ for which
$\omega(X)$ is constant are not $H$-invariant.  However, it is not
dif\/f\/icult to show that $\Omega$ is constant in Sharpe's sense if and
only $DT=0$.  Let $\{e_i\}$ be a~basis for $\g$, $\{E_i\}$ the vector
f\/ields on $P$ such that $\omega(E_i)=e_i$.  Then
$\Omega(E_i,E_j)=\Omega^k_{ij}e_k$, and in a~space of constant
curvature the coef\/f\/icients $\Omega^k_{ij}$ are constants.  Let
$C^k_{ij}$ be the structure constants of~$\g$ relative to the basis
$\{e_i\}$.  Then $C^k_{ij}-\Omega^k_{ij}$ are the structure constants
for a new Lie algebra stucture on (the vector space) $\g$, making it
into a new Lie algebra $\g'$.  This process of changing the Lie
algebra bracket is called mutation.  Note that since $\Omega$ vanishes
if one of its arguments is vertical, mutation does not change $\h$.

Starting with a Cartan space form modelled on $(\g,\h)$ with group
$H$, by carrying out the procedure in Theorem~\ref{theorem3} one obtains, for the
same Cartan algebroid, the Lie-algebraic bracket appropriate to a
Cartan geometry modelled on $(\g',\h)$ with group $H$, where $\g'$ is
the mutation of $\g$ determined by the (constant) curvature $\Omega$.
With respect to the new Lie-algebraic bracket we have
$T(\sigma,\tau)=-\{\sigma,\tau\}$, so the mutant geometry has
curvature zero.  This illustrates an observation of Sharpe's, that all
Cartan space forms are mutants of f\/lat geometries.  Blaom puts the
same idea in a dif\/ferent way:\ his theory, he says, `{\em has no
models\/}:\ all symmetric structures are created equal and the
curvature [of $\nablaC$] \ldots\ merely measures deviation from {\em
some\/} symmetric structure' (his italics). The tractor connection
$\nablaT$, on the other hand, does dif\/ferentiate between dif\/ferent
Cartan space forms.

\subsection*{Acknowledgements}
I am grateful to one of the referees for several helpful comments,
and in particular for drawing my attention to reference \cite{Cap2}
and pointing out its relevance.

I am a Guest Professor at Ghent University:\ I should like to express
my gratitude to the members of the Department of Mathematical Physics
and Astronomy for their hospitality.

\pdfbookmark[1]{References}{ref}

\LastPageEnding


\begin{thebibliography}{99}

\footnotesize\itemsep=0pt

\bibitem{Blaom}
Blaom A.D.,
Geometric structures as deformed inf\/initesimal symmetries,
{\it Trans.\ Amer.\ Math.\ Soc.} {\bf 358} (2006), 3651--3671,
\href{http://arxiv.org/abs/math.DG/0404313}{math.DG/0404313}.

\bibitem{Cap2}
\v{C}ap A.,
Inf\/initesimal automorphisms and deformations of parabolic geometries,
{\it J. Eur. Math. Soc. (JEMS)} {\bf 10} (2008), 415--437,
\href{http://arxiv.org/abs/math.DG/0508535}{math.DG/0508535}.

\bibitem{Cap}
\v{C}ap A., Gover A.R.,
Tractor calculi for parabolic geometries,
{\it Trans.\ Amer.\ Math.\ Soc.} {\bf 354} (2001), 1511--1548.

\bibitem{Sharpe}
Sharpe R.W.,
Dif\/ferential geometry. Cartan's generalization of Klein's Erlangen program, with a foreword by S.S.~Chern, {\it Graduate Texts in Mathematics}, Vol.~166, Springer-Verlag, New York, 1997.


\end{thebibliography}
\end{document}